\documentclass[11pt]{article}

\pdfoutput=1
\usepackage{graphicx}
\usepackage{amsmath,amssymb,amsfonts}%
\usepackage{amsthm}
\usepackage[colorlinks=true,
linkcolor=black,
citecolor=blue,
urlcolor=black]{hyperref}
\usepackage{xcolor}

\newtheorem{proposition}{Proposition} 

\newtheorem*{definition}{Definition}%

\title{ 
Addition to "Structured random matrices and cyclic cumulants: A free probability approach"}

\author{Denis Bernard\thanks{Laboratoire de Physique de Ecole Normale Sup\'erieure, CNRS, ENS \& PSL University, Sorbonne Universit\'e, Universit\'e Paris Cit\'e, 24 rue Lhomond, Paris, 75005, France}, \
Ludwig Hruza\thanks{Laboratoire de Neurosciences Cognitives et Computationnelles, INSERM U960, Ecole Normale Superieure-PSL Research University, Paris, France}\ .}

\begin{document}

\maketitle

Abstract: 
We give a refined definition of the class of random matrix ensembles introduced in our paper \cite{bernard_StructuredRandom_2024} by extending the so-called fourth axiom to deal with cumulants of disjoint cycles. We argue that the theorems concerning the stability of such ensembles under non-linear transformations still hold with these refined axioms.

\vskip 0.5 truecm

In our paper \cite{bernard_StructuredRandom_2024} we had defined a class of structured random matrix ensembles by postulating four axioms, and showed the stability of such matrix ensembles under a set of non-linear transformations. Here we replace axiom (iv) by a stronger version which makes a statement about cumulants of disjoint cycles, which in particular allows to deal with cumulants of products of traces of powers of the matrix. This new axiom is for example satisfied by the ensemble of symmetric orthogonally \cite{speicher_EntrywiseApplication_2024} and unitarily invariant random matrices \cite{Collins_cumulants_2007}. For QSSEP we have good reasons to believe that axiom (iv) is correct, but we have not proved it yet.
The refined version of the axioms is:

\begin{definition}[Axioms of structured random matrix ensembles]
	\label{axioms}
We consider ensembles of random hermitian matrices $M=M^\dag$ with measure $\mathbb E$ that satisfy, in the large $N$ limit, the four following defining properties~:
\begin{enumerate}
	\item[(i)] Local $U(1)$-invariance, meaning that, in distribution, $M_{ij}\stackrel{d}{=}e^{-i\theta_i}M_{ij}e^{i\theta_j}$, for any phases $\theta_i$, $\theta_j$.
	\item[(ii)] Cyclic cumulants of $n$ matrix elements scale as $N^{1-n}$, meaning that\\ $C_n[M_{i_1 i_2},M_{i_2i_3},\cdots ,M_{i_n i_1}]=\mathcal O(N^{1-n})$, for any choice of $i_k$.
	\item[(iii)] Scaled cyclic cumulants  $N^{n-1}C_n[M_{i_1 i_2},M_{i_2i_3},\cdots ,M_{i_n i_1}]$ are continuous functions in $x_k=i_k/N$ at coinciding indices, in the large $N$ limit.
	\item[(iv)] Cumulants of $r$ disjoint cycles with $n_1+\cdots+n_r =: n$ elements in total scale as $N^{2-r-n}$. That is,
	\[
	C_n[\underbrace{M_{i_1^{1} i_2^{1}},\cdots, M_{i_{n_1}^{1}i_1^{1}}}_{\text{$1$st cycle with $n_1$ elements}},\cdots , \underbrace{M_{i_1^{r} i_2^{r}},\cdots, M_{i_{n_r}^{r}i_1^{r}}}_{\text{$r$th cycle with $n_r$ elements}}] \sim N^{2-r-n}.
	\]
	as long as no pair of cycles has an index in common.
\end{enumerate}
\end{definition}

Since these axioms are formulated in terms of cumulants, they are preserved by shifting the matrices $M$ by diagonal matrices, $M_{ij}\to M_{ij}-\delta_{ij} a_i$. Without loss of generality we can assume the matrices to be centered, $\mathbb{E}[M_{ii}]=0$. 

As in \cite{bernard_StructuredRandom_2024}, let us introduce the so-called local free cumulants $g_n(x_1,\cdots,x_n)$ defined as the scaled cyclic cumulants,
\begin{equation}
g_n(x_1,\cdots,x_n) := \lim_{N\to\infty} N^{n-1}C_n[M_{i_1 i_2},M_{i_2i_3},\cdots ,M_{i_n i_1}] ~,
\end{equation}
with $x_k=i_k/N$ fixed in this limit. The dependence of $g_n$'s on the points $x_k$ reflects the structure of these random matrices.

Let us recall that the first three axioms (i-iii) imply a simple formula for the expectation values of single cyclic traces. For any diagonal matrix $\Delta_k$ with entries $(\Delta_k)_{ii}=\Delta_k(i/N)$ with $\Delta_k(\cdot)$ smooth functions, we write by linearity in the $\Delta$'s, 
\begin{equation}\label{eq:expectation_trace}
\lim_{N\to\infty} \mathbb{E}[\underline{\mathrm{tr}}(M\Delta_1M\Delta_2\cdots M\Delta_n)] = \int_0^1 \prod_{k=1}^n dx_k \Delta_k(x_k)\, T_n(x_1,\cdots,x_n) ~,
\end{equation}
with $\underline{\mathrm{tr}}(\cdot):=N^{-1}\mathrm{tr}(\cdot)$ the normalized trace. 
Then  \cite{bernard_StructuredRandom_2024}, the functions $T_n$'s  can be written as sum over the set $NC(n)$ of non-crossing partitions of order $n$ as
\begin{equation}
T_n(x_1,\cdots,x_n) = \sum_{\pi \in NC(n)} g_\pi(x)\, \delta_{\pi^*}(x) ~,
\end{equation}
where $g_\pi(x)=\prod_{b\in\pi} g_b(x)$ the local free cumulants associated to $\pi\in NC(n)$, $\pi^*$ is the Kreweras dual of $\pi$ and $\delta_{\pi^*}$ the delta function forcing the points in each block of $\pi^*$ to be equal.

The new version of axioms (iv) in particular implies that 
\begin{equation}
\mathbb{E}[e^{zN\,\underline{\mathrm{tr}}M}]\asymp_{N\to\infty} e^{zN\,\mathbb{E}[\underline{\mathrm{tr}} M]}~,
\end{equation}
as a formal power in $z$, so that the trace of $M$ is self-averaging. This is shown by expanding $\log\mathbb{E}[e^{z\,N\,\underline{\mathrm{tr}}M}]$ in terms of cumulants $\frac{z^nN^n}{n!}C_n[\underline{\mathrm{tr}}M,\cdots,\underline{\mathrm{tr}}M]$ and evaluating the scaling at each order as a function of $N$ using axiom (iv).

The main aim of this note is to show that these axioms, including the new version of axiom (iv), are preserved by non-linear transformations.

\begin{proposition}[Invariance under polynomial operations]\label{Prop1}
	The axioms, in particular axiom (iv), are invariant under matrix-valued polynomial operations $M \to P(M)$, for any polynomial $P$.
\end{proposition}

\begin{proposition}[Invariance under entry-wise operations]\label{Prop2}
	Assume that the matrices are centered, $\mathbb{E}[M_{ii}]=0$. Then, the axioms, in particular axiom (iv), are invariant under entry-wise polynomial operations $M_{ij} \to Y_{ij}= M_{ij}\,f_{ij}(N|M_{ij}|^2)$, for any non-linear power series $f_{ij}(\cdot)$ that can be different for each entry, with a finite large $N$ limit $f_{xy}(\cdot)$.
\end{proposition}

Before going to the proof, let us note that, as corollary of the first proposition, the trace of any polynomial in $M$ is self-averaging, that is
\begin{equation}
\mathbb{E}[e^{zN\,\underline{\mathrm{tr}}P(M)}]\asymp_{N\to\infty} e^{zN\,\mathbb{E}[\underline{\mathrm{tr}}P(M)]}~,
\end{equation}
for any polynomial $P$. This has numerous applications. 

Another corollary is that
	\begin{equation} \label{eq:trace-cumulants}
		C_r[\underline{\mathrm{tr}}(M^{n_1}),\cdots,\underline{\mathrm{tr}}(M^{n_r})] \sim N^{2-2r}.
	\end{equation}
	which can be obtained from the first proposition using the polynomial  $P(M)=a_1M^{n_1}+\cdots+a_rM^{n_r}.$ 
But contrary to the case $r=1$ in Eq.~\eqref{eq:expectation_trace}, we don't have formulas for the finite coefficient on the right hand side of Eq.~\eqref{eq:trace-cumulants}, since this would need more information than what is contained in the local free cumulants $g_n$.

\begin{proof}
	Here we only consider axiom (iv), for the others we refer to \cite{bernard_StructuredRandom_2024}. As will be come clear in the end of "Justification of Proposition 1", we are lacking a formal proof of an important claim, that is based on insights we gained from doing many examples. 
	
	The argument is going to be based on using three operations on partitions such that the partitions contributing non-trivially to the cumulants of $r$ cycles can be induced by recursive applications of these operations from simple partitions for which we control the scaling in $N$. These operations can be viewed as moves on the lattice of partitions. We first apply this method to justify the first Proposition and then the second one.
	
	\subparagraph{Graphical representation.}
	Let us consider $r$ disjoint cycles consisting of a total of $n$ matrix elements. Each matrix element is represented by a line between two black dots (nodes are indices and edges are matrix elements),
	\begin{equation}
		\raisebox{-0.5\height}{\includegraphics{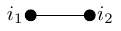}} := M_{i_1i_2}.
	\end{equation} 
	Then replace any matrix element $M_{i_1i_2}$ by an arbitrary power $(M^p)_{i_1i_2}$ and represent this by an intermediate white circle\footnote{One could also draw a white circle for each additional matrix element and sum over the additional indices. Then, instead of one circle, there would be $p$ circles here. But this complicates things, because every of these white circles behaves the same.}
	\begin{equation}\label{eq:segment}
		\raisebox{-0.5\height}{\includegraphics{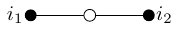}} := \sum_{j_1,\cdots,j_p} M_{i_1j_1}M_{j_1j_2}\cdots M_{j_pi_2}.
	\end{equation} 
	Now assume that in each of the $r$ disjoint cycles, one replaces each of the \mbox{$n_1 + \cdots + n_r =: n$} matrix elements by an arbitrary power, such that in total there are now $m_1+\cdots +m_r =: m$ additional matrix elements, together with $m$ sums over the newly introduced indices.
	To prove the first proposition with this representation (suppressing the indices on the black dots, which here are required to be distinct across cycles), we would like to show that 
	\begin{equation}\label{eq:disjoint_cycles}
		C_n\left[\; \raisebox{-0.5\height}{\includegraphics{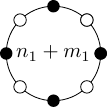}}\quad\cdots\quad \raisebox{-0.5\height}{\includegraphics{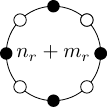}}\;\right]\sim N^{2-r-n},
	\end{equation}
	where the commas in the cumulant $C_n$ are located at the black dots, such that the cumulant is of order $n$, even though there are $n+m$ matrix elements in total.
	
	\subparagraph{Cumulant formula.}To reduce such a cumulant of products of variables to cumulants of the individual variables, we recall a formula from \cite{leonov_MethodCalculation_1959}. Denote $k=k_{1}+\cdots+k_{n}$ and $\Gamma=\gamma_{1}\cup\cdots\cup\gamma_{n}\in P(k)$ a partition of $\{1,\cdots,k\}$ into intervals $\gamma_{1},\cdots,\gamma_{n}$ of length $k_{1},\cdots,k_{n}$. For $\pi\in P(k)$, denote $\pi\vee\Gamma$ the maximum of the two partitions, then
	\begin{equation} \label{eq:cumulants-product}
		C_{n}[X_{1}^{k_{1}},\cdots,X_{n}^{k_{n}}]=\sum_{\pi\in P(k):\,\pi\vee\Gamma=1_{k}} \!\! C_{\pi}[\underbrace{X_{1},\cdots,X_{1}}_{k_{1}\text{ times}},\cdots,\underbrace{X_{n},\cdots,X_{n}}_{k_{n}\text{ times}}] ~.
	\end{equation}
	In other words, one has to sum over all partitions $\pi$ (crossing and non-crossing) whose parts $p\in\pi$ successively connect the intervals $\gamma_{1},\cdots,\gamma_{n}$ to become the complete interval $1_{k}$. 	
	
	\subparagraph{Justification of Proposition \ref{Prop1}.}
	In our case, an interval $\gamma $ corresponds to a line segment between two black dots with an intermediate white circle, such as in Eq.~\eqref{eq:segment}. There are $n$ such segments in total, and the total number of matrix elements is $k=n+m$.	
		
	The strategy is to start from a factorized partition that does not yet satisfy $\pi\vee \Gamma = 1_{n+m}$, but for which we know the scaling with $N$, and to successively modify this partition until it satisfies the condition.
	
	We start with $\pi=\pi_1\cup\cdots\cup\pi_r$ which is the union of disjoint partitions $\pi_c$ of the individual cycles $c=1,\cdots,r$, each respecting $\pi_c\vee\Gamma|_c = 1_{n_c+m_c}$ within the respective cycle $c$. From \cite{hruza_CoherentFluctuations_2023}, we know that in this case the leading contribution of $C_\pi$ comes from partitions $\pi$ where each $\pi_c$ is non-crossing. Furthermore, the Kreweras complements $\pi_c^*$ determine which indices (these can be either black or white dots, but only within the same cycle) must be equal to respect local $U(1)$ invariance. As a result, all $|p_c|$ elements in a part  $p_c\in\pi_c$ are arranged into a single loop and therefore the part contributes with $N^{1-|p_c|}$. Taking into account the constraints on the $m_c$ sums in this cycle enforced by $\pi_c^*$, the cumulant scales as $C_{\pi_c}\sim N^{1-n_c}$. The overall scaling of $C_\pi$ is then $N^{1-n_1}\times\cdots\times N^{1-n_r}=N^{r-n}$. 
	This scaling has to be compared to the cost induced by the minimum number of operations that are necessary to connect all $r$ cycles, according to Eq.~\eqref{eq:cumulants-product}. If each operation reduces the number of disjoint cycles by one, then the minimum number of necessary operations is $r-1$. If each of these operations "cost" a factor of $N^{-2}$, the resulting partition will scale as $N^{r-n}\times N^{-2(r-1)}=N^{2-r-n}$ exactly as required.
	If one of the operations costs less (a higher negative power of $N$), for example $N^{-3}$, the partition will be subleading. If an operation costs more, for example $N^{-1}$, this would falsify our proposition. 
	
	To avoid confusion, we will use the word \textit{cycle} only for the $r$ big cycles in Eq.~\eqref{eq:disjoint_cycles}, and we use the word \textit{loop} for cycles that appear within any part of $\pi$. By local $U(1)$ invariance, any part of $\pi$ must close into a loop to give a non-zero contribution. The operations we are discussing below preserve these conditions.
	
	\textbf{(1)} The first possible operation is to join two arbitrary parts $p_{c_1}\in\pi_{c_1}$ and $p_{c_2}\in\pi_{c_2}$ from two separate cycles into a single part $p_{c_1}\cup p_{c_2}$. Originally, the two separate parts contributed with a scaling of $N^{1-|p_{c_1}|}N^{1-|p_{c_2}|}$. After joining the parts, the new part consists of two loops and contributes $N^{-|p_{c_1}|-|p_{c_2}|}$, by axiom (iv). This results in an additional scaling factor of $N^{-2}$.
	Schematically, and abbreviating $[\cdots]\equiv C_n[\cdots]$,
	\begin{equation}\label{eq:operation-1}
		\left[\; \raisebox{-0.5\height}{\includegraphics{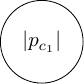}}\right]
		\left[\; \raisebox{-0.5\height}{\includegraphics{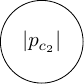}}\right]
		\overset{N^{-2}}{\longrightarrow}
		\left[\; \raisebox{-0.5\height}{\includegraphics{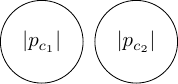}}\right].
	\end{equation}
	where the simple circles stand for an arbitrary number of matrix elements $M_{ij}$ connected in a loop, without explicitly showing black or white dots.
	
	\textbf{(2)} For the second operation, one selects two parts $p_{c_1}\in\pi_{c_1}$ and $p_{c_2}\in\pi_{c_2}$ from distinct cycles, divides these parts into segments $b_{1}\subset p_{c_1}$ and $b_{2}\subset p_{c_2}$ and their complements $\bar b_{1}$ and $\bar b_{2}$, and defines two new parts as $b_{1}\cup b_{2}$ and $\bar b_{1}\cup \bar b_{2}$. Schematically, without explicitly showing black or white dots, except on the endpoints of the segments where the gray dot can be either black or white (but at least one of the two dots connected by a Kronecker delta must be white), 
	\begin{equation}\label{eq:join-segments}
		\left[\; \raisebox{-0.5\height}{\includegraphics{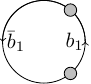}}\right]
		\left[\; \raisebox{-0.5\height}{\includegraphics{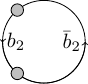}}\right]
		\overset{N^{-2}}{\longrightarrow}
		\left[\; \raisebox{-0.5\height}{\includegraphics{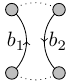}}\right]
		\left[\; \raisebox{-0.5\height}{\includegraphics{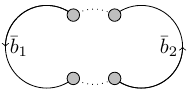}}\right].
	\end{equation}
	 Here local $U(1)$ invariance requires two new Kronecker deltas, represented by dotted lines. This operation therefore "costs" a factor of $N^{-2}$. 
	 
	 If we had only joined $b_{1}\cup b_{2}$ into a new part, but had kept  $\bar b_{1}$ and $\bar b_{2}$ as separate parts, then this case can be reduced to the first operation, by starting with partitions $\pi_{c_1}$ and $\pi_{c_2}$ whose parts $p_{c_1}$ and $p_{c_2}$ have already been divided into subparts, together with the necessary Kroneker deltas that ensure local $U(1)$ invariance of the subparts within each cycle separately.
	 
	 However, if one would select disjoint segments in the same cycles, such as $b_1,b_1'$ in
	 \begin{equation}
	 	\left[\; \raisebox{-0.5\height}{\includegraphics{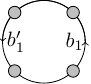}}\right]
	 \end{equation}
	 to form a new part $b_1\cup b_1'\cup b_2$, then more Kronecker deltas are necessary and one can check that the operation would cost less (a higher negative power) than $N^{-2}$.
	 Also note that if one applies the second operation to the special case where $b_1=p_1$, $b_2=p_2$ and empty complements, this becomes,
	 	 \begin{equation}
	 	\left[\; \raisebox{-0.5\height}{\includegraphics{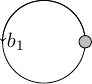}}\right]
	 	\left[\; \raisebox{-0.5\height}{\includegraphics{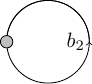}}\right]
	 	\overset{N^{-2}}{\longrightarrow}
	 	\left[\; \raisebox{-0.5\height}{\includegraphics{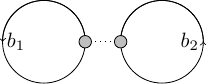}}\right]
	 \end{equation}
	 which contributes to the final result, but in a different form than the outcome of operation 1 in Eq.~\eqref{eq:operation-1}.
	 
	\textbf{(3)} To get all possible allowed partitions by successive iteration of the moves, one needs a third operation which combines parts from three different cycles as follows (for clarity we have inserted the indices $i_1,j_1,i_2,\cdots$ at the boundary of the segments)
	\begin{align}\label{eq:operation-3}
		\left[\; \raisebox{-0.5\height}{\includegraphics{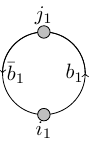}}\right]
		\left[\; \raisebox{-0.5\height}{\includegraphics{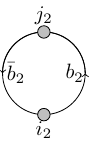}}\right]
		\left[\; \raisebox{-0.5\height}{\includegraphics{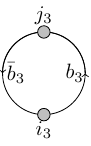}}\right]
		\overset{N^{-4}}{\longrightarrow}
		\left[\; \raisebox{-0.5\height}{\includegraphics{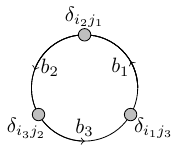}}\right]
		\left[\; \raisebox{-0.5\height}{\includegraphics{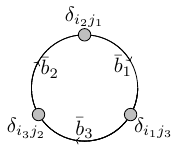}}\right].
	\end{align}
	Directly connecting three cycles together (doing two steps in one) the operation costs $N^{-4}$. This is because the left hand side scales as $N^{3-k_\mathrm{tot}}$ (with $k_\mathrm{tot}$ the total number of matrix in these three parts), the right hand side scales as $N^{2-k_\mathrm{tot}}$, and one has to introduce $3$ Kronecker deltas to ensure $j_1=i_2$, $j_2=i_3$ and $j_3=i_1$. Note that the operation reduces to operation (2) in case the third part $b_3 \cup \bar b_3$ is empty.
	
	Operation (3) can also be understood as a concatenation of operation (2) and an operation (3') that opens up the previously created Kronecker delta $\delta_{i_1j_2}$ and inserts the third part by creating two new deltas $\delta_{j_2i_3}$ and $\delta_{j_3i_1}$. From this point of view, the intermediate step is in Eq.~\eqref{eq:operation-3} is
	
	\begin{align}
		\cdots \quad
		\underset{(2)}{\overset{N^{-2}}{\longrightarrow}}
		\left[\; \raisebox{-0.5\height}{\includegraphics{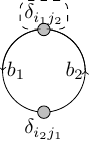}}\right]
		\left[\; \raisebox{-0.5\height}{\includegraphics{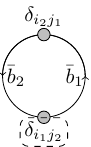}}\right]
		\left[\; \raisebox{-0.5\height}{\includegraphics{tikz/op3b_part3}}\right]
		\underset{(3')}{\overset{N^{-2}}{\longrightarrow}}
		\cdots 
		\quad
	\end{align}


To conclude, we claim that that operations 1, 2 and 3 are the only relevant operations at leading order, and that there is no operation that is less costly. Having considered many examples, we have good reasons to believe this. However, we are lacking an actual proof, which would probably involve a more formal argument about the lattice of partitions and the set of partitions we can reach with our three operations.

\subparagraph{Justification of Proposition \ref{Prop2}.}
This proof is actually simpler as there are no indices to sum over. Recall the notation $Y_{ij}=M_{ij}\,f_{ij}(N|M_{ij}|^2)$. By linearity, we may consider $f_{ij}(x)=x^{m_{ij}}$ with arbitrary degrees $m_{ij}$ for each edge. Cumulants w.r.t $Y$ of disjoints cycles are computed in terms of cumulants of $M$ using Eq.~\eqref{eq:cumulants-product}. Here the number of $M$'s in a given cycle $c$ is $k_c=n_c+2m_c$ with $m_c:=m_{i_1^{c}i_2^{c}}+\cdots m_{i_{n_c}^{c}i_1^{c}}$, where $i_l^{c}$, with $l=1,\cdots,n_c$, denotes the indices of the cycle $c$.

Again, we start from a partition $\pi=\pi_1\cup\cdots\cup\pi_r$, which is the union of disjoint partitions $\pi_c$ of the individual cycles $c=1,\cdots,r$, each respecting $\pi_c\vee\Gamma|_c = 1_{k_c}$. For each cycle $c$ with indices $i_l^{c}$, $l=1,\cdots,n_c$, the leading contribution to the cumulant $C_{\pi_c}$ of $M$ comes from forming a long loop $M_{i_1^{c}i_2^{c}}\cdots M_{i_{n_c}^{c}i_1^{c}}$, and a collection of small loops of length two of the form $M_{i_l^{c}i_{l+1}^{c}}M_{i_{l+1}^{c}i_l^{c}}$. There are $m_c$ such small loops. Each long loop scales as $N^{1-n_c}$ and each small loop as $\mathcal{O}(1/N)$. The total scaling $N^{-m_c}$ of all small loops gets compensated by the factors $N^{m_{ij}}$ arising from the definition  $Y_{ij}=M_{ij}\,f_{ij}(N|M_{ij}|^2)$] so that the global scaling of $\pi$ is $N^{r-n}$.

We then apply iteratively the moves (1-2-3) to $\pi$ to enforce the condition $\pi \vee\Gamma= 1_k$, where $k=k_1+\cdots+k_r$. Since there are no indices to sum over and indices on distinct cycles are different, only the moves of type (1) matters. There are various ways to apply those moves, depending whether we group long or small loops together. But for any of these ways we gain a factor $N^{-2}$. Again $(r-1)$ operations of this type are needed to enforce $\pi \vee\Gamma= 1_k$, so that at the end the global scaling is $N^{r-n}N^{-2(r-1)}=N^{2-r-n}$ as required.
\end{proof}

We end this note by listing  a few misprints observed in our manuscript \cite{bernard_StructuredRandom_2024} (these misprint have been corrected in the latest arXiv version v5):
\begin{itemize}
	\item In the paragraph above Eq.~(4): \\
"with the trace of the dressed matrix $M_h=h^{\frac{1}{2}}Mh^{\frac{1}{2}}$'' should be replaced by "with the matrix element $v_h^TMv_h$, with $v_h$ the vector with coordinates $h_i^\frac{1}{2}$'' -- since $\mathrm{tr}(QM)=v_h^TMv_h$ for $q(x,y)=h^{\frac{1}{2}}(x)h^{\frac{1}{2}}(y)$ --, and ''$F_0[h]:=\lim_{N\to\infty}N^{-1}\log\mathbb{E}[e^{N\underline{tr}(M_h)}]$'' should be replaced by ''$F_0[h]:=\lim_{N\to\infty}N^{-1}\log\mathbb{E}[e^{v_h^TMv_h}]$.
	\item In Proposition 3:\\
As formulated in the text, this Proposition holds assuming the matrices $M$ to be centered: $\mathbb{E}[M_{ii}]=0$.
This is not restrictive, as the axioms are invariant under shifting $M$ by diagonal matrices, but it avoids having to define specifically the diagonal matrix elements $Y_{ii}$.
\end{itemize}

\bibliography{ref-comment_axiom_4}
\bibliographystyle{plain}

\end{document}